\documentclass[12pt]{article}
\usepackage{amsmath}
\usepackage{amsthm}
\usepackage{amssymb}
\usepackage{amscd}
\usepackage[mathscr]{euscript}

\begin{document}
\newtheorem{Lemma}{Lemma}[section]
\newtheorem{Theorem}[Lemma]{Theorem}
\newtheorem{Proposition}[Lemma]{Proposition}
\newtheorem{Corollary}[Lemma]{Corollary}
\newtheorem{Definition}[Lemma]{Definition}
\newtheorem{Example}[Lemma]{Example}
\newtheorem{Remark}[Lemma]{Remark}

\newcommand{\bthe}{\begin{Theorem}}
\newcommand{\ethe}{\end{Theorem}}
\newcommand{\bpr}{\begin{Proposition}}
\newcommand{\epr}{\end{Proposition}}
\newcommand{\bco}{\begin{Corollary}}
\newcommand{\eco}{\end{Corollary}}
\newcommand{\bre}{\begin{Remark} \em}
\newcommand{\ere}{\end{Remark}}
\newcommand{\ble}{\begin{Lemma}}
\newcommand{\ele}{\end{Lemma}}
\newcommand{\bex}{\begin{Example} \em}
\newcommand{\eex}{\end{Example}}
\newcommand{\bde}{\begin{Definition} \em}
\newcommand{\ede}{\end{Definition}}
\newcommand{\bpf}{\begin{proof}[\bf{Proof.}]}
\newcommand{\epf}{\end{proof}}
\newcommand{\bqu}{\begin{Question} \em}
\newcommand{\equ}{\end{Question}}

\renewcommand{\th}{\theta}
\renewcommand{\t}{\tilde}
\newcommand{\e}{\eta}
\newcommand{\Om}{\Omega}
\newcommand{\om}{\omega}
\newcommand{\al}{\alpha}
\newcommand{\si}{\sigma}
\newcommand{\De}{\Delta}
\newcommand{\de}{\delta}
\renewcommand{\b}{\beta}
\newcommand{\ga}{\gamma}
\newcommand{\Ga}{\Gamma}
\newcommand{\ep}{\varepsilon}
\newcommand{\La}{\Lambda}
\newcommand{\la}{\lambda}
\newcommand{\va}{\varphi}

\newcommand{\rg}{\rangle}
\renewcommand{\lg}{\langle}
\newcommand{\cc}{\mathds{C}}
\newcommand{\rr}{\mathbb{R}}
\newcommand{\qq}{\mathbb{Q}}
\newcommand{\zz}{\mathbb{Z}}
\newcommand{\nn}{\mathbb{N}}
\newcommand{\ff}{\mathbb{F}}
\newcommand{\End}{{\rm End}}
\newcommand{\Hom}{{\rm Hom}}
\newcommand{\id}{{\rm id}}
\newcommand{\rad}{{\rm rad}}
\newcommand{\Rad}{{\rm Rad}}
\newcommand{\Mod}{{\rm Mod}}
\newcommand{\im}{{\rm im}}
\newcommand{\Ext}{{\rm Ext}}
\newcommand{\soc}{{\rm soc}}
\newcommand{\Soc}{{\rm Soc}}
\newcommand{\ann}{{\rm ann}}
\newcommand{\Spec}{{\rm Spec}}
\newcommand{\mxa}{\mathcal{M}(X,\mathcal{A})}
\newcommand{\nxa}{\mathcal{M}'(X,\mathcal{A})}
\newcommand{\lra}{\longrightarrow}
\newcommand{\ma}{\mathcal{A}}
\newcommand{\mb}{\mathcal{B}}
\newcommand{\mxb}{\mathcal{M}(X,\mathcal{B})}
\newcommand{\mxc}{\mathcal{M}(X,\mathcal{C})}
\newcommand{\mra}{\mathcal{M}(\rr,\mathcal{A})}
\newcommand{\mna}{\mathcal{M}(\nn,\mathcal{A})}
\newcommand{\mc}{\mathcal{C}}
\newcommand{\mf}{\mathcal{F}}
\newcommand{\mxbs}{\mathcal{M}(X,\mathcal{B}_S)}
\newcommand{\m}{{\rm m}}
\newcommand{\mrb}{\mathcal{M}(\rr,\mathcal{B})}
\newcommand{\mrl}{\mathcal{M}(\rr,\mathcal{L})}
\newcommand{\mxz}{\mathcal{M}(X,Z(X))}
\newcommand{\mxcz}{\mathcal{M}(X,Coz(X))}

\title{Remarks On $\aleph_0$-Injectivity}

\author{E. Momtahan \\
\texttt{e-momtahan@mail.yu.ac.ir}}

\date{\small{Department of Mathematics, Yasouj University, Yasouj, Iran}}

\maketitle

\begin{abstract}
\noindent
 We continue our review of the literature via countable injectivity.
\end{abstract}

\vspace{0.5cm}


\vspace{0.5cm}
\section{Introduction}
In this article, by $R$ we mean an associative ring with identity
and $M$ is  a unitary $R$-module. An $R$-module $M$ is said to be
$\aleph_0$-injective (f-injective) if every (module) homorphism
$f\in \rm{Hom}_{R}(I,M)$, there exists $\bar{f}\in
\rm{Hom}_{R}(R,M)$, such that $\bar{f}|_{I} = f$, where $I$ is
any countably generated (finitely generated) right ideal of $R$. A
ring $R$ is said to be right $\aleph_0$-self-injective
(f-injective) if it is $\aleph_0$-injective (f-injective) as a
right $R$-module. In this note, by $C(X)$ we mean the ring of all
 real valued continuous over a completely regular space (or
 equivalently a Tychonoff space). The reader is referred to
 \cite{gj} for undefined notations and definitions. A module is
 called extending ($\aleph_0$-extending) if every submodule
 (countably generated submodule) is a direct summand.

 This paper is a continuation of \cite{mom1} and
\cite{mom2}, in which the author studied $\aleph_0$-injectivity
of modules and rings. In \cite{mom1}, many well-known
observations on injective modules have been seen  to be true if
we replace injectivity by $\aleph_0$-injectivity. Among them, for
example, one can quote, this one:  a ring $R$ is semisimple
artinian if and only every $R$-module (or every countably
generated $R$-module) is $\aleph_0$-injective. Using this result,
one can prove the following proposition which is well-known if we
replace injectivity by $\aleph_0$-injectivity. The equality of
the first two parts of the following result is due to B. Osofsky
(see \cite[page 114, excersise 11]{la}).

\bpr The following assertions are equivalent:
\begin{enumerate}

\item $R$ is semisimple Artinian

\item the intersection of any two injective $R$-modules is injective.

\item the intersection of any two $\aleph_0$-injective $R$-modules is $\aleph_0$-injective.

\end{enumerate}

\epr

\bpf The proof which is given for $(1)\Leftrightarrow(3)$ works
also for $(1)\Leftrightarrow(2)$. $(1)\Rightarrow (3)$ is obvious.
Suppose that the intersection of any two $\aleph_0$-injective
$R$-modules is $\aleph_0$-injective. Now suppose that $R$ were
not semisimple Artinian.  In \cite{mom1}, it has been observed
that a ring is semi-simple Artinian if and only if every module is
$\aleph_0$-injective. Hence by our assumption there is an
$R$-module $M$ which is not $\aleph_0$-injective. Set
$E=E(M)\oplus E(E(M)/M)$. Then $E$ is injective. Let $x$ belong to
$E(M)\setminus M$. Set $y=(x, x+M)\in E$. Then $(M,0)$ is
essential in $(E(M),0)$ and $E((M,0)+yR)\subseteq E$ and the
intersection of those two injectives is $M$ which is not
$\aleph_0$-injective. A contradiction.

\epf

\bre It is well-known that a ring is von Neumann regular if and
only if every module over $R$ is f-injective. Now the same method
used in the above result  shows that (i) a ring $R$ is von Neumann
regular if and only if the intersection of any two f-injective
module is f-injective. It is also well-known that a ring $R$ is
(left and right) artinian serial with $Jac(R)^2=0$ if and only if
every module over $R$ is extending (see \cite[13.5]{du}). Again,
the same line of proof shows that a ring $R$ is (left and right)
artinian serial if and only if the intersection of any two
extending modules is extending. \ere

However, if every finitely generated (or even cyclic) $R$-module
is $\aleph_0$-injective, we do not need to access the
semisimplicity of $R$ (contrary to Osofsky's famous result which
asserts that a ring $R$ is semisimple artinian if and only if
every cyclic $R$-module is injective). In fact, a ring is right
$\aleph_0$-self-injective regular if and only if every cyclic
$R$-module is $\aleph_0$-injective.  In \cite{mom2}, it has been
found out that the well-known result of Y. Utumi that {\it every
right (or left) self-injective ring module its
 Jacobson radical is von Neumann regular}, is no longer true if we
 replace injectivity by $\aleph_0$-injectivity.
 A commutative $\aleph_0$-self-injective ring with many additional
 properties has been constructed that is not von Neumann regular modulo its
 Jacobson radical.

 The next proposition will help us in the sequel. It is
well-known for injectivity (see \cite{wis1}).

\bpr Let $M$ be an $R$-module and $I\subseteq \rm{Ann}(M)$. If
$M$ is an $\aleph_0$-injective $R$-module, then $M$ is an
$\aleph_0$-injective $R/I$-module. The converse is true provided
that $R$ is a right fully idempotent ring. \epr

\bpf Let $M$ be an $\aleph_0$-injective $R$-module, and
$f:A/I\longrightarrow M$ be an $R/I$-homomorphism, where $A/I$ is
a countably generated right ideal. Thus we may put $A=B+I$, where
$B=\sum_{k=1}^{\infty}a_kR$ and then let $f(a_k+I)=m_k$, where
$m_k\in M$. Defining $g:B\longrightarrow M$ by
$g(a_k)=m_k=f(a_k+I)$, by our hypothesis $g$ can be extended to
$R$. Suppose now $M$ is $\aleph_0$-injective as an $R/I$-module,
and let $f:A\longrightarrow M$ be an $R$-homomorphism, where $A$
is a countably generated right ideal. First we claim that
$f(A\bigcap I)=0$: for every $x\in A\bigcap I$ we have
$x=\sum_{k=1}^{k=n}a_kb_k$ where $a_i$ and $b_i$ belong to
$A\bigcap I$ and therefore
$$f(x)=f(\sum_{k=1}^{k= n}a_kb_k)=\sum_{k=1}^{k=n}f(a_k)b_k$$
but $f(a_k)b_k \in MI=0$ so $f(A \bigcap I)=0$. In as much as
$A+I/I \cong A/A \bigcap I$ and $f(A \bigcap I)=0$ we infer that
$\bar{f}:A+I/I\longrightarrow M$ is an $R/I$-module homomorphism
so $\bar{f}:R/I\longrightarrow M$, i.e.,
$\bar{f}(a+I)=(a+I)(r+I)=ar+I$ for some $r\in R$, i.e. $f(a)=ar$.
\epf

\section{$\aleph_0$-Ikeda-Nakayama rings}
In \cite[Theorem 1]{ik}, it has been proved that a ring $R$ is
right f-injective if and only if  (i)$l(I_1\bigcap
I_2)=l(I_1)+l(I_2)$, where $I_1$ and $I_2$ are finitely generated
left ideals of $R$ and (ii) $r(l(I))$ where $I$ is a principal
ideal right ideal of $R$. Rings satisfy in (i), for each pair of
right ideals are called Ikeda-Nakayama (IK-)rings. Such rings have
been extensively studied (see for example \cite{cam}). IK-rings
have been shown to be a special class of quasi-continuous rings
(see \cite{cam}, where it has been observed for the first time).
Furthermore, a theory of IK-modules is shown in \cite{wis2002},.
It is easy to see that every right self-injective module is a
right IK-ring. However, the converse is not true. This is
because,  the ring of integers, $\Bbb{Z}$, is an IK-ring that is
not self-injective, since $\Bbb{Z}$ is not a divisible abelian
group. Here we consider $\aleph_0$-IK rings and study their
behaviour in  Boolean rings and rings of real valued
 continuous functions over a Tychonoff space.

\bde
 A ring is called a right $\aleph_0$-IK ring if $l(I_1\bigcap
I_2)=l(I_1)+l(I_2)$, where $I_1$ and $I_2$ are countably
generated left ideals of $R$. \ede

We begin with a lemma which is useful in the sequel:

\ble If $R$ is a commutative $\aleph_0$-selfinjective ring, then
$R$ is an $\aleph_0$-IK ring.
 \ele

\bpf Obviously $l(I_1)+l(I_2)\subseteq l(I_1\bigcap I_2)$.
Suppose on the other hand that $a\in l(I_2\bigcap I_2)$, we can
define a homomorphism $\alpha: I_1+I_2\longrightarrow R$ as:

\[ \alpha(b) = \left \lbrace
 \begin{array}{c l}
    b & b\in I_1,\\
    (1+a)b & b\in I_2.
  \end{array}
\right. \]

Since these two expressions coincide on $I_1\bigcap I_2$, by
$\aleph_0$-injectivity of $R$, $\exists c\in R$ such that
$\alpha(b)=cb$ for $b\in I_1$. Thus, we have $cb=b$, i.e.,
$(c-1)b=0$. Consequently we  write $a=(c-1)+(1+a-c)$($c-1\in
l(I_1)$ and $(1+a-c)\in l(I_2)$), this proves the lemma.
 \epf

 For the sequel we
need the following result which is  an slight modification of
Theorem 8 in \cite{wis2002}.

\ble Let $M_R$ be a right $R$-module and $S=\rm{End}(M)$. Then
the following are equivalent:
\begin{enumerate}
\item $M$ is $\aleph_0$-$\pi$-injective ($\aleph_0$-quasi
continuous)
\item For any two countably generated submodules $A$ and $B$ of
$M_R$, $S=l_{S}(A)+l_{S}(B)$.

\end{enumerate}
\ele

Using these results we have the following result in $C(X)$. In
the following a space $X$ is called {\it extremally disconnected}
if open sets have open closure. By a {\it basically disconnected}
space we mean  a space in which co-zero sets (i.e. the complement
of zero sets) have open closure.

\bthe Let $X$ be a completely regular space. Then the following
are equivalent:
\begin{enumerate}
\item for any two countably generated ideals $A$ and $B$ of
$C(X)$,with $A\bigcap B =(0)$;
$$\rm{Ann}(A)+\rm{Ann}(B)=C(X)$$
\item $X$ is basically disconnected.
\end{enumerate}

\ethe

\bpf (1)$\Rightarrow$ (2) By Lemma 1.2, we observe that $C(X)$ is
$\aleph_0$-continuous. Hence it is an $\aleph_0$-extending ring.
By \cite{az}, Theorem 3.3, $X$ is basically disconnected.

(2)$\Rightarrow$ (1) Since $X$ is basically disconnected, again
by \cite{az}, Theorem 3.3, we observe that $C(X)$ is
$\aleph_0$-extending. Let $\langle e_1 \rangle$ and $\langle e_2
\rangle$ be two ideals which are summands of a commutative ring
$R$. It is well-known that  $\langle e_1 \rangle + \langle e_2
\rangle = \langle e_1+e_2-e_1e_2\rangle$ and $\langle e_1\rangle
\bigcap \langle e_2\rangle=\langle e_1 e_2 \rangle$. Hence $C(X)$
satisfies summand sum property (SSP) and summand intersection
property. This implies that $C(X)$ is an $\aleph_0$-extending
ring if and only if $C(X)$ is an $\aleph_0$-quasi continuous ring.

\epf

\bthe Let $X$ be a completely regular space. Then the following
are equivalent:
\begin{enumerate}
\item for any two ideals $A$ and $B$ of $C(X)$ with $A\bigcap B
=(0)$,
$$\rm{Ann}(A)+\rm{Ann}(B)=C(X)$$
\item $X$ is an extremally disconnected space,
\end{enumerate}
\ethe

\bpf (1)$\Rightarrow$(2) By \cite{wis2002}, Theorem 8,  $C(X)$ is
a quasi-continuous ring and hence extending. By \cite{az}, Theorem
3.5, we deduce that $X$ is extremally disconnected.

(2)$\Rightarrow$(1): By \cite{az}, Theorem 3.5. when $X$ is
extremally disconnected, then $C(X)$ is an extending ring. But as
we stated in the proof of the above theorem, $C(X)$ has always
summand sum property. Therefore it is quasi-continuous when it is
extending. Now by \cite[Theorem 8]{wis2002} the implication
follows. This completes the proof.
 \epf

\bre It is natural to speculate on the space $X$ if $C(X)$
satisfies  the following stronger conditions: (i) for any two
ideals $A$ and $B$ of $C(X)$, $\rm{Ann}(A)+\rm{Ann}(B)=C(X)$, and
(ii) for any two countably generated ideals $A$ and $B$ of $C(X)$,
$\rm{Ann}(A)+\rm{Ann}(B)=C(X)$. We conjecture that in the first
case $X$ is an extremally $P$-space (and hence $C(X)$ is
self-injective regular in this case) and in the second case, $X$
would be a $P$-space. We leave these questions open for those who
are interested in the theory of rings of continuous functions. In
\cite{es}, it has been shown that $C(X)$ is
$\aleph_0$-self-injective if and only if $C(X)$ is regular. Hence
if our conjecture (i.e., conjecture (ii)) is true, we would have:
$C(X)$ is $\aleph_0$-IK if and only if $C(X)$ is
$\aleph_0$-self-injective.
 \ere

 The following result is due to O. A. S. Karamzadeh and A. A.
Koochakpour \cite{ka}. Let $A$ and $B$ be two subsets of the ring
$R$. We say $A\bigcup B$ is an {\it orthogonal set} if $\forall
x,y\in A\bigcup B$,  $x\neq y$, then $xy=0$. We say that an
element $x$ separates $A$ from $B$ if $xa^2=a$ for all $a\in A$
and $xB=0$. If there exists such an $x$, then we say that $A$ has
a left separation from $B$.

\ble For a strongly regular ring $R$ the following are equivalent:
\begin{enumerate}
\item $R$ is $\aleph_0$-self-injective

\item If $S\bigcup T$ is a countable orthogonal set in $R$ with
$S\bigcap T= \emptyset$, then $S$ has a left separation from $T$.

\end{enumerate}
\ele

Andrew B. Carson \cite{ca}, has shown that if $R$ is an
$\aleph_0$-complete Boolean ring (Boolean joins of countably many
elements always exist), then $R$ is $\aleph_0$-self-injective. He
also showed that there are $\aleph_0$-self-injective Boolean
rings which are not $\aleph_0$-complete. In spite of this fact,
we show that in a Boolean ring the following fact holds:

\bpr Let $R$ be a Boolean ring, then the following are
equivalent:
\begin{enumerate}

\item $R$ is an $\aleph_0$-IN-ring.

\item for any two ideals $A$ and $B$ of $R$ with $A\bigcap B
=(0)$, $\rm{Ann}(A)+\rm{Ann}(B)=R$.

\item $R$ is $\aleph_0$-self-injective.

\end{enumerate}
\epr

 \bpf
  (1)$\Rightarrow$(2) It is always true.\

(2)$\Rightarrow$(3): We want to show that $R$ is
$\aleph_0$-self-injective ring. According to Lemma 2.7, it is
enough to show that each two disjoint orthogonal countable
subsets of $R$ can be separated by an element of $R$. Let
$S\bigcup T$ be an orthogonal set with $S\bigcap T=\emptyset$.
Since $R$ is an $\aleph_0$-IN ring and $Ann(B)=Ann(\langle B
\rangle)$, where $B\subseteq R$ and by $\langle B \rangle$ we
mean the ideal generated by the set $B$, we have $$Ann( \langle S
\rangle \bigcap \langle T \rangle)=Ann(\langle S
\rangle)+Ann(\langle T \rangle)$$

It is then evident that $R=Ann(\langle S \rangle\bigcap \langle T
\rangle)=Ann(\langle S \rangle)+Ann(\langle T \rangle)$, i.e.
$1=x+y$, where $xS=0$ and $yT=0$. Now we see that $1-x=y$
separates $S$ from $T$, because $yT=0$ and $yS=(1-x)S=S$ i.e.
$ya=(1-x)a=a-xa=a$. But $R$ is a Boolean ring, hence $ya^2=a$.\

(3)$\Rightarrow$(1) follows from  Lemma 2.1. \epf

Along this line, we observe that a Boolean ring is  an IK-ring if,
and only if, it is self-injective. We need the following lemma
which has been proved in \cite{karam}:

\ble For a strongly regular ring $R$ the following are equivalent:
\begin{enumerate}
\item $R$ is self-injective

\item If $S\bigcup T$ is an orthogonal set in $R$ with
$S\bigcap T= \emptyset$, then $S$ has a left separation from $T$.

\end{enumerate}
\ele

\bpr Let $R$ be a Boolean ring, then the following are equivalent:
\begin{enumerate}

\item $R$ is an IN-ring,

\item for any two ideals $A$ and $B$ of $R$ with $A\bigcap B
=(0)$, $\rm{Ann}(A)+\rm{Ann}(B) = R$.

\item $R$ is self-injective.

\end{enumerate}
\epr

\bpf Use Lemma 2.9 and follow exactly the proof of Proposition
2.8. \epf
\section{$\sum$-$\aleph_0$-injectivity}

In \cite{fa}, C. Faith has shown that a $\sum$-injective
$R$-module $M$ with endomorphism ring $S$ is characterized by the
ascending chain condition on the lattice of $S$-submodules which
are annihilators of subsets of $R$ (\cite[Prposition 3.3]{fa}).
If $E(R)$ denotes the injective hull of $R_R$, and if $M=E(R)$,
this condition implies the ascending chain condition on
annihilator right ideals (=right annulets) of $R$, and in case
$M=E(R)=R$, this condition is equivalent to a.c.c. on right
annulets (\cite[Corollary 3.4 and Theorem 3.5]{fa}). Now we
observe that some parts of these results by C. Faith in \cite{fa}
can be applied to answer a similar question concerning
$\sum$-$\aleph_0$-injectivity. Let $M_R$ be a left $R$-module
 and $S=End(M)$. Following  Faith's nomenclature, for any subset $X$ of
 $M$,

 $$X^{\bot}=\rm{Ann}_{R}(X)=\{r\in R\;|\;Xr=0\}$$\\
 is a right ideal of $R$. The set of such right ideals is denoted
 by $\ma$. By $\ma_{\aleph_0}$, we mean all countably
 generated members of $\ma$. And by $\mb_{\aleph_0}$, we mean the
 set of all annihilators of members of $\ma_{\aleph_0}$ in $M$.
 Since $I\longrightarrow I^{\bot}$ is an injective and order-inverting
 map between $\ma_{\aleph_0}$ and $\mb_{\aleph_0}$,  one satisfies
 the ascending chain condition if, and only if, the other one satisfies the descending chain
 condition.

 \bpr
If $M_R$ is an $R$-module, then $\ma_{\aleph_0}$ satisfies the
ascending chain condition if and only if for each countably
generated right ideal $I$ of $R$ there exists a finitely
generated subideal $I_1$ such that $I^{\bot}=I_1^{\bot}$.
\epr

\bpf Assume a.c.c. for $\ma_{\aleph_0}$ or equivalently, the
d.c.c. for $\mb_{\aleph_0}$. Let $I$ be a countably generated
right ideal of $R$, and let $I_1$ be a finitely generated
subideal of $I$, such that $I_1^{\bot}$ be a minimal element of
the set
$$\{ K^{\bot}: \mbox{$K$ is a finitely generated subideal
of $I$} \}$$ in $M$. If $x\in I$, then $I_1+xR$ is also a finitely
generated subideal of $I$, and hence satisfies
$(I_1+xR)^{\bot}\subseteq I_1^{\bot}$. By the choice of $I_1$,
necessarily $(I_1+xR)^{\bot}=I_1^{\bot}$, so $I_1^{\bot}x=0$.
Since this is true for all $x\in I$, then $I_1^{\bot}I=0$,
whence, $I_1^{\bot}\subseteq I^{\bot}$. Moreover, $I_1\subseteq I$
implies $I^{\bot}\subseteq I_1^{\bot}$, consequently
$I_1^{\bot}=I^{\bot}$ follows.

Conversely, let $I_1\subseteq I_2\subseteq \cdots \subseteq I_n
\subseteq \cdots$ be a chain of countably generated right ideals
of $R$ lying in $\ma_{\aleph_0}$, let $X_i=I_i^{\bot}$,
$i=1,2,\cdots$, be the corresponding elements of
$\ma_{\aleph_0}$, and suppose also $I=\bigcup
_{n=1}^{\infty}I_n$. Now, let $J$ be the finitely generated
subideal if $I$ such that $I^{\bot}= J^{\bot}$. Since $J$ is
finitely generated, there is an integer $q$ such that $J\subseteq
I_k$, $k\geq q$. Thus, $I_k^{\bot}\subseteq J^{\bot}$, $k\geq q$.
Moreover,
$$J^{\bot}=I^{\bot}=\bigcap_{n=1}^{\infty}I_n^{\bot},$$

Consequently, $I_k^{\bot}=J^{\bot}$, for $k\geq q$. Then,
$I_k=(I_k^{\bot})^{\bot}=I_q$, $k\geq q$, and the result follows.

 \epf

\bco A ring $R$ satisfies the a.c.c. on $\mb_{\aleph_0}$ if and
only if each  countably generated right ideal $I$ contains a
finitely generated ideal $I_1$, such that $I^{\bot}=I_1^{\bot}$.
\eco

\bpr The following conditions on an $\aleph_0$-injective module
$M$ are equivalent:
\begin{enumerate}

\item $M^{(\nn)}$ is $\aleph_0$-injective,

\item $R$ satisfies the a.c.c. on the ideals in $\ma_{\aleph_0}$,

\item $M$ is $\sum$-$\aleph_0$-injective.

\end{enumerate}
 \epr

 \bpf
(1)$\Rightarrow$(2) Let $I_1\subseteq I_2\subseteq
\cdots\subseteq I_m\subseteq \cdots$ be a strictly ascending
chain of right ideals in $\ma_{\aleph_0}$, let
$I=\bigcup_{n=1}^{\infty} I_n$, and let $x_n$, be an element of
$I_n^{\bot}$ (taken in $M$) not in $I_{n+1}^{\bot}$,
$n=1,2,\cdots$. If $r\in I$, then $x_nr=0$ for all $k\geq q$.
Therefore the element $r'=(x_1r,\cdots, x_nr, \cdots)$ lies in
$M^{(\nn)}$, even though $x=(x_1,\cdots,x_n,\cdots)$ lies in
$M^{\nn}$. Let $f$ denote a map defined by $f(r)=r'$ for all
$r\in I$. Assuming that $M^{(\nn)}$ is injective, there exists,
by Baer's criterion, an element $y=(y_1,\cdots,y_m,0,\cdots)\in
M^{(\rr)}$ such that
$$f(r)=yr=(y_1r,\cdots,y_mr,0,\cdots)=(x_1r,\cdots, x_mr,\cdots)$$

for all $r\in I$. But this implies that $x_tr=0$ for all $t>m$,
for all $r\in I$, that is, $x_t\in I^{\bot}\subseteq
I_{t+1}^{\bot}$, this contradicts the choice of $x_t$.

(2)$\Rightarrow$(3). Let $I$ be a right ideal of $R$, and
$I_1=r_1R+\cdots+r_nR$ be the finitely generated subideal given
by the above proposition such that $I^{\bot}=I_1^{\bot}$. Let
$f:I\longrightarrow M^{(\Lambda)}$ be any map. Since $M^{\Lambda}$
is injective (it is direct product of injective modules), there
exists an element $p\in M^{\Lambda}$ such that $f(r)=pr$ for all
$r\in I$. Since $f(r_i)=pr_i\in M^{(\Lambda)}$, $i=1,\cdots, n$,
there exists an element $p'\in M^{(\Lambda)}$ such that
$p_ar_i=p'_ar_i$ for all $a\in \Lambda$, $i=1,\cdots,n$, where
$g_a$ is the $a$ coordinate of any $g\in M^{\Lambda}$. Since
$r_1,\cdots,r_n$ generate $I_1$, this implies that $pr=p'r$ for
all $r\in I_1$, hence $(p_a-p'_a)\in I_1^{\bot}$ for all $a\in
\Lambda$. Since $I_1^{\bot}=I^{\bot}$, it follows that for all
$x\in I$ and for all $a\in A$, that is $px=p'x$, for all $x\in
I$. Thus, $f(x)=p'x$ for all $x\in I$, with $p'\in
M^{(\Lambda)}$,  so $M^{(\Lambda)}$ is $\aleph_0$-injective by
Baer's criterion. \epf

Along this line and as an application of Proposition 1.3, we
consider the following theorem which is well-known when we
replace $\aleph_0$-injectivity by injectivity. We need  the
following lemma.

\ble Suppose $R$ is a regular ring and $M$ is an $R$-module. Then
$M$ is  $\sum$-$\aleph_0$-injective if and only if $R/A$ is an
artinian ring where $A=\rm{Ann}(M)$. \ele

\bpf A standard proof given in \cite{pa}, works here with only
 a slight modification. \epf

 \bthe Let $R$ be a
regular ring which is an algebra over a field $F$. Let $M$ be an
$R$-module with $\rm{dim}_{F} M \leq \aleph_0$, and let
$I=\rm{Ann}(M)$. Then the following are equivalent:
\begin{enumerate}

\item $R/I$ is Artinian

\item $M$ is $\sum$-$\aleph_0$-injective

\item $M$ is $\aleph_0$-injective.

\end{enumerate}

 \ethe

 \bpf
Using the  proof of the above theorem (see  \cite{pa}) and apply
Proposition 1.3.
 \epf

\section{On $\aleph_0$-quasi injectivity}

It is well-known that $\aleph_0$-injectivity is a Morita property
for regular rings (\cite[Chapter 14]{go}). In fact, by a theorem
of D. Handelman we have: Let $R$ be a right
$\aleph_0$-self-injective regular ring. If $M$ is a finitely
generated projective right $R$-module, then $T=End(M)$ is a right
$\aleph_0$-self-injective regular ring. We say that a module is
{\it $\aleph_0$-quasi injective} if every $R$-homomorphism $f:
B\longrightarrow M$ extends to $M$, when $B$ is a countably
generated submodule of $M$. By the proof of Handelman's theorem
we infer that:

\bco Let $R$ be an $\aleph_0$-selfinjective regular ring, and $P$
is a finitely generated projective module, then $P$ is
$\aleph_0$-quasi injective.
 \eco

As we saw over $\aleph_0$-self-injective regular rings, every
finitely generated free module is $\aleph_0$-quasi injective.
This is a good motivation to study $\aleph_0$-quasi injectivity.
J. Ahsan \cite{ah}  introduced the concept of  qc-ring. Here we
will present results on $\aleph_0$-qc ring.

\bde A ring $R$ is called right $\aleph_0$-qc ring if every
$R$-homomorphic image of $R$ as a right $R$-module is
$\aleph_0$-quasi injective.
\ede

The next result is an extension of a theorem of Nicholson and
Yousif \cite{ni}. Although it is stated for quasi injective
module, with only a slight modification, it is true for
$\aleph_0$-quasi injective modules. Recall that a module $M$ is
Dedekind-finite if $M\oplus N\cong M$ then $N=0$.

\bpr Let $M$ be a quasi injective module. Then $M$ is
Dedekind-finite if, and only if, $M$ is co-hopfian. \epr

\bpf It is well-known that $M$ is Dedekind-finite if, and only
if, $\rm{End}(M)$ is a Dedekind-finite ring. Now let
$f:M\longrightarrow M$ be an $R$-homomorphism. If $f$ is a
monomorphism, then there exits $g: f(M)\longrightarrow M$ such
that $g=f^{-1}$, but $f(M)$ is a submodule of $M$, so there
exists $g':M\longrightarrow M$ such that $g'f=1$ (by quasi
injectivity of $M$). But $\rm{End}(M)$ is Dedekind-finite so
$fg'=1$, whence, $f$ is an epimorphism. Conversely, let
$fg=1_{M}$. This implies that $g$ is monic ($g(m)=0\Rightarrow fg
(m)=1(m)=m=0$). By $M$ is co-hofian, i.e., $g$ is an epimorphism,
i.e., $\exists h$ such that $gh=hg=1$. It follows that $h=f$, and
$gf=1$.
 \epf

\ble Let $R$ be a ring, and $I$ an ideal of $R$. If $M$ is an
$\aleph_0$-quasi injective $R/I$-module, then $M$ is also
$\aleph_0$-quasi injective as an $R$-module. Also, if $M$ is an
$\aleph_0$-quasi injective $R$-module, then $M$ is also
$\aleph_0$-quasi injective as an $R/I$-module \ele

\bpf The relation $(*)$: $m(r+I)=mr$($m\in M$ and $r\in R$) is
used in each case to define $M$ as a module over $R$ or $R/I$,
where $M$ is given as a module over $R/I$ or $R$. Consider that
if $N=\langle x_1, x_2, \cdots \rangle$ is countably generated as
a an $R/I$-module, then for all $x\in N$ there exist
$i_1,i_2,\cdots, i_k$ such that $x=x_{i_1}(r_{i_1}+I)+\cdots +
x_{i_n}(r_{i_n}+ I)$. By $(*)$, $x=x_{i_1}r_{i_1}+\cdots+
x_{i_n}r_{i_n}$, this means that $N$ is also countably generated
as an $R$-module. It is easy to see that the concepts "submodule"
and "homomorphism" coincide over each ring. Hence any diagram

$$\begin{array}{ccccc}
  0 & \longrightarrow & A & \longrightarrow & M \\
   & & \downarrow &  &  \\
   &  & M &  &
\end{array}$$\\
over any ring is also a diagram over the other ring. Thus $M$ is
$\aleph_0$-quasi injective over $R$ if, and only if, $M$ is
$\aleph_0$-quasi injective over $R/I$ module. The result follows.

\epf

\ble Let $R$ be a ring and $K$ an arbitrary two sided ideal of
$R$. Then $R$ is a right $\aleph_0$-qc ring if, and only if,
$R/K$ is a right $\aleph_0$-qc ring. \ele

\bpf Suppose $R$ is a right $\aleph_0$-qc ring and $K$ is a two
sided ideal of $R$. We show that $R/K$ is a right $\aleph_0$-qc
ring. Let $I/K$ be any right ideal of $R/K$ with $I\subseteq R$
and consider the right $R/K$-module $R/K/I/K\cong R/I$ as an
$R$-module with $K\subseteq I$. Then $K$  annihilates the
$R$-module $R/I$, and therefore, $R/I$ may be regarded as an
$R/K$-module. Furthermore, $R$ is a right $\aleph_0$-qc ring by
hypothesis. Hence, $R/I$ is $R$-$\aleph_0$-quasi-injective. Hence,
by the above lemma, $R/I$, as an $R/K$-module, is
$R/K$-$\aleph_0$-injective. We have shown that any cyclic
$R/K$-module is $R/K$-$\aleph_0$-quasi-injective. Therefore $R/K$
is a right $\aleph$-qc ring.\epf

\bthe Let $R$ be a commutative ring. Then $R$ is an $\aleph_0$-qc
ring if, and only if, every factor ring of $R$ is an
$\aleph_0$-self-injective ring. \ethe

\bpf If $I$ is an ideal of an $\aleph_0$-qc ring $R$, then $R/I$
is an $\aleph_0$-qc ring by lemma 3.5. Therefore $R/I$ is an
$\aleph_0$-self-injective ring (it is evident that since $R_R$ is
generated by the identity, we may infer that any homomorphism from
countably generated ideal of $R$ into $R_R$ can be extended to an
endomorphism). Hence $R_R$ is $\aleph_0$-injective, and $R$ is a
cyclic $R$-module. Then, $M\cong R/I$ for some ideal $I$ of $R$.
By assumption, $R/I$ is $R/I$-$\aleph_0$-quasi-injective. Hence
$R/I$ is $R$-$\aleph_0$-quasi-injective, and therefore $R$ is an
$\aleph_0$-qc ring. \epf

\section{Examples}
 In this section we provide some new examples of
 $\aleph_0$-self-injective (regular) rings which are perhaps of some interest for
 their own right.

\bex Let $F$ be a field and $G$ a group. If $\rm{char}(F)=0$, then
the following are equivalent:
\begin{enumerate}

\item $F[G]$ is $\aleph_0$-self-injective;

\item $G$ is finite;

\item $F[G]$ is self-injective.

\end{enumerate}
\eex

\bpf (1)$\Rightarrow$(2): Suppose $F[G]$ is
$\aleph_0$-self-injective, then it is $p$-injective. By a result
of Farkas (\cite{pa}), $G$ is locally finite. Now by a result of
Villamayor-Connell \cite{pa}, page 69, Theorem 1.5, $F[G]$ is a
regular ring. The rest of proof is the same as Renault's Theorem
(\cite{pa}, Theorem 2.8). The remaining implications are
well-known and can be found in  \cite{pa}. \epf

\bex Let $X$ be a Tychonoff space. In \cite{es}, it has been shown
that $C(X)$ is $\aleph_0$-self-injective if, and only if,
$C(X)/C_{F}(X)$ is $\aleph_0$-self-injective if, and only if, $X$
is a P-space and therefore, if, and only if, $C(X)$ is regular.
\eex

\bex Let $\ma$ be a $\si$-Algebra. In \cite{azad}, it has been
observed that rings of all real valued $\ma$-measurable functions
are $\aleph_0$-self-injective. \eex

\bex Let $X$ be Tychonoff space. By $D(X)$, we mean the lattice of
continuous
    functions $f$ on $X$ with values in the extended real numbers
    $\rr\bigcup\{\pm \infty\}$, for which $f^{-1}\rr$ is a dense
    subset of $X$. In general, under pointwise addition and
    multiplication, $D(X)$ is not a ring. However, when $X$ is a
    quasi-F space, then $D(X)$ is a ring. A Tychonoff space $X$ is called a quasi F-space
    if every dense co-zero set $S$ in $X$ is $C^{\ast}$-embedded.   Since $D(X)\cong D(\beta
    X)$, we may without loss of generality suppose that $X$ is a
    compact space. Then by a result due to Hager (see \cite{azad}),
     $D(X)\cong \frac{\mxa}{N}$, for some certain $N$. This then implies
    that $D(X)$ is also an $\aleph_0$-self-injective regular ring, for
     $X$ a quasi-F space (see \cite{azad})
\eex

The following example is a slight modification of a theorem by F.
L. Sandomierski (see \cite{sa}).

\bex Let $M$ be a right $R$-module and a direct sum of countably
many non-zero submodules $\{M_n\;|\;n\in \nn\}$ and
$S=\rm{End}(M)$. Then $_{S}M$ is not an $\aleph_0$-injective
$S$-module. \eex

\bpf Let $e_i:M_R\longrightarrow M_i$ be the $i$-th projection of
the module $M_R$ onto the submodule $M_i$, then $\{e_i\}_{i\in
\nn}$ is a countable set of orthogonal idempotents of $S$. Let
$_{S}A$ be the ideal of $S$ generated by $\{e_i\}_{i\in \nn}$.
Since $M_i\neq 0$ for each $i\in I$, choose $0\neq x_i\in M_i$.
Clearly there is an $S$-homomorphism $f:_{S}A\longrightarrow
_{S}A$ such that $e_if=x_i=e_ix_i$. If $f$ were extendable to a
homomorphism from $S$ to $_{S}M$, then it would be given by some
element of $M$. However, for any element $x\in M$, $e_ix=0$ for
all but finitely many $i\in \nn$, so $f$ is not extendable to
$_{S}S$, so $_{S}M$ is not $\aleph_0$-injective. \epf

\end{document}